\documentstyle[amscd]{amsart}
\newtheorem{thm}{Theorem}[section]
\newtheorem{lm}[thm]{Lemma}
\newtheorem{cor}[thm]{Corollary}
\newtheorem{pro}[thm]{Proposition}

\theoremstyle{definition}
\newtheorem{df}[thm]{Definition}

\theoremstyle{remark}
\numberwithin{equation}{section}

\def \R {\bold{R}}
\def \Z {\bold{Z}}
\def \CA {\cal A}
\def \CB {\cal B}
\def \CC {\cal C}
\def \CG {\cal G}

\def \CL {\cal L}
\def \CM {\cal M}
\def \hM {\hat{{\cal M}}}
\def \CR {\cal R}
\def \CP {\cal P}

\def \la {\langle}
\def \ra {\rangle}
\def \a {\alpha}

\def \c {\chi}
\def \g {\gamma}
\def \G {\Gamma}
\def \d {\delta}

\def \Lam {\Lambda}
\def \lam {\lambda}

\def \ov {\overline}
\def \op {\oplus}
\def \p {\phi}

\def \s {\sigma}
\def \Si {\Sigma}
\def \t {\tau}
\def \T {\Theta}

\def \O {\Omega}

\def \bd {\partial}

\def \x {\times}
\def \ox {\otimes}
\def \ve {\varepsilon}
\def \D {\Delta}

\def \e {\eta}

\begin{document}

\baselineskip.525cm

\title[Monopole Homology]
{A monopole homology for integral homology 3-spheres}

\author[Weiping Li]{Weiping Li}
\address{Department of Mathematics, Oklahoma State University \newline
\hspace*{.175in}Stillwater, Oklahoma 74078-0613}
\email{wli@@math.okstate.edu}
\date{February 8, 2000}

\begin{abstract}
To an integral homology 3-sphere $Y$, we assign a well-defined
$\Z$-graded (monopole) homology 
$MH_*(Y, I_{\e}(\T; \e_0))$ whose construction in principle follows from the 
instanton Floer theory with the dependence of the spectral flow
$I_{\e}(\T; \e_0)$, where $\T$ is the unique $U(1)$-reducible monopole of 
the Seiberg-Witten equation on $Y$ and $\e_0$ is a reference perturbation
datum. 
The definition uses the moduli space of
monopoles on $Y \x \R$ introduced by Seiberg-Witten in studying smooth
4-manifolds. We show that the monopole homology $MH_*(Y, I_{\e}(\T; \e_0))$ 
is invariant
among Riemannian metrics with same $I_{\e}(\T; \e_0)$. This provides a
chamber-like structure for the monopole homology of
integral homology 3-spheres. The assigned function
$MH_{SWF}: \{I_{\e}(\T; \e_0)\} \to \{MH_*(Y, I_{\e}(\T; \e_0))\}$ is a
topological invariant (as Seiberg-Witten-Floer Theory).
\end{abstract}

\maketitle

\section{Introduction}

Since Donaldson \cite{dk} initiated the study of smooth
4-manifolds via the Yang-Mills theory, the gauge theory (Donaldson
invariants, relative Donaldson-Floer invariants and Taubes'
Casson-invariant interpretation, etc) has proved remarkably fruitful
and rich to unfold some of the mysteries in studying smooth 4-manifolds.
The topological quantum field theory proposed by Witten \cite{wit}
stimulates the
most exciting developments in low-dimensional
topology. In 1994, Seiberg and Witten introduces a new (simpler) kind of
differential-geometric equation (see \cite{don, wi}). 
In a very short time after the equation
was introduced, some long-standing problems were solved, new and unexpected
results were discovered. For instance, Kronheimer and Mrowka \cite{km}
proved the Thom conjecture affirmatively, several authors
proved variants (generalizations) of the Thom conjecture
independently in \cite{fs, mst, ru}, as well as the three-dimensional
version of the Thom conjecture \cite{au}. Taubes showed that
there are more constraints on symplectic structures in \cite{tau, tau1}
and the beautiful equality $SW = Gr$ in \cite{tau2, tau3}.
See \cite{don} for a survey in the Seiberg-Witten theory.

Using the dimension-reduction principle, one expects the 
Floer-type homology of 3-manifolds via the Seiberg-Witten equation.
Indeed Kronheimer and Mrowka \cite{km} analyzed the 
Seiberg-Witten-Floer theory for
$\Si \x S^1$, where $\Si$ is a closed oriented surface.
Later on Marcolli studied the Seiberg-Witten-Floer homology for
3-manifolds with first Betti number positive in \cite{m}.
For a connected compact oriented 3-manifold with positive first Betti
number and zero Euler characteristic, Meng and Taubes \cite{mt} showed that
a (average) version of Seiberg-Witten invariant is the same as the Milnor
torsion.
The interesting class of 3-manifolds as integral (rational) homology 3-spheres
is lack of well-posed theory.
Although various authors attempted to resolve the problem on
defining a ``Seiberg-Witten-Floer'' theory, 
the new phenomenon of harmonic-spinor
jumps and the dependence of Riemannian metrics is not addressed clearly.
The metric-dependence (also related to the harmonic-spinors) issue is
quickly realized by many experts in this field (see \cite{don, moy}). In
\cite{moy}, the irreducible Seiberg-Witten-Floer homology of Seifert
space is shown to be dependent on the metric and the choice of connection
on the tangent bundle (as our reference $\e_0$ in this paper).

In this paper, we construct a monopole homology from the Seiberg-Witten
equation in the same way as an instanton Floer homology from
the Self-Duality equation in Donaldson-Floer theory \cite{fl}.
Our key point is that by using the
unique $U(1)$-reducible solution $\T$ of the Seiberg-Witten equation on
an integral homology 3-sphere $Y$ we make use of the spectral flow of $\T$ 
to capture the dependence in certain perturbation classes of 
Riemannian metrics and
1-forms. The same idea was used before by the present author to
establish a symplectic Floer homology of knots in \cite{li}, and the
original one was in the study of the instanton Floer homology of
rational homology 3-spheres by Lee and the present author in \cite{ll}.
Many technique issues such as transversality, transitivity and gluing property
are treated in many authors books and papers, those techniques follow
the same line in \cite{fl} or simpler. So we omit the details on these, but
only emphasize the Riemann-metric dependence and understand the role of
such a fixing spectral flow of $(\T; \e_0)$.

Our approach is similar to approaches in \cite{clm, ll, li} to understand
the perturbation data (including Riemannian metrics).
The unique $U(1)$-reducible $\T$ gives a spectral flow 
$I_{\e}(\T; \e_0)$ as a Maslov index in \cite{clm} Part III. The spectral flow
$I_{\e}(\T; \e_0) = \mu_{\e}(\T) - \mu_{\e_0}(\T)$ with respect to a reference
$\e_0$ fixes a class of admissible perturbations consisting of
Riemannian metrics and 1-forms. As long as Riemannian metrics and 1-forms
give the same spectral flow $I_{\e}(\T; \e_0)$, we
prove that the constructed monopole homology is invariant inside
the fixed class of Riemann-metrics and 1-forms $(\e = (g_Y, \a))$ with same 
$I_{\e}(\T; \e_0)$.
The spectral flow $I_{\e}(\T; \e_0)$
is not a topological invariant, and is dependent upon the
Riemannian metrics. Without fixing a class of Riemannian metrics with
same $I_{\e}(\T; \e_0)$, one cannot obtain well-defined notions such as
spectral flow of irreducible Seiberg-Witten solutions on $Y$, and the gluing
formula as well as the relative Seiberg-Witten invariant. Hence our results
follow from 
fixing $I_{\e}(\T; \e_0)$.

\noindent{\bf Theorem A. }{\em
(1) For an integral homology 3-sphere $Y$ and any admissible perturbation $\e$,
there is a well-defined $\Z$-graded monopole homology
$MH_*(Y, I_{\e}(\T; \e_0))$ constructed 
by the Seiberg-Witten equation over $Y\x {\R}$.

(2) For any two admissible perturbations $\e_1$ and $\e_2$, there is a
group homomorphism $\Psi_*$ between two monopole
homologies $MH_*(Y, I_{\e_1}(\T; \e_0))$ and $MH_*(Y, I_{\e_2}(\T; \e_0))$. 

(3) If $I_{\e_1}(\T; \e_0) = I_{\e_2}(\T; \e_0)$, then the homomorphism
$\Psi_*$ is an isomorphism.}

Our fixed-class $I_{\e}(\T; \e_0)$ of Riemannian metrics gains control of the
birth and death of irreducible solutions of the Seiberg-Witten equation on the 
integral homology 3-sphere $Y$. Changing the reference $\e_0$ into $\e_0^{'}$
corresponds to the overall degree-shifting by $\mu_{\e_0^{'}}(\T) - 
\mu_{\e_0}(\T)$ for monopole homologies.
The control in the instanton homology of 
rational homology 3-spheres is gained by fixing the spectral flows of 
all $U(1)$-reducibles from the Wilson-loop perturbations (not metrics).
The control in the monopole homology of integral homology 3-spheres
is gained by fixing the spectral flow of the unique $U(1)$-reducible
$\T$ from the Riemannian metrics (not only 1-forms).
Fixing $I_{\e}(\T; \e_0)$ enters crucially in proving Theorem A and 
Theorem B.

\noindent{\bf Theorem B. }{\em
For a smooth 4-manifold $X =X_0 \#_Y X_1$
with $b_2^+(X_i) > 0 (i = 0, 1)$ and $Y$ an integral homology 3-sphere,
the Seiberg-Witten invariant of $X$ is
given by the Kronecker pairing of $MH_*(Y; I_{\e}(\T; \e_0))$ with
$MH_{-1-*}(-Y;I_{\e}(\T; \e_0))$ for the relative Seiberg-Witten invariants
$q_{X_0,Y, \e}$ and $q_{X_1,- Y, \e}$ (see Definition~\ref{relative});
\[ \langle  , \rangle :
MH_*(Y;I_{\e}(\T; \e_0)) \x MH_{-1-*}(-Y;I_{\e}(\T; \e_0)) \to {\bf Z}; 
\ \ \  q_{SW}(X) =
\langle q_{X_0,Y, \e}, q_{X_1,-Y, \e}\rangle. \]}

The paper is organized as follows. \S 2 provides an introduction of
the Seiberg-Witten equation on 3-manifolds. \S 3 studies the
configuration space over $Y$ through Seiberg-Witten equation and a
natural monopole complex. We show that there are admissible
perturbations from 
Riemannian metrics and 1-forms in \S 4 via the method similar to
\cite{pt}. The spectral-flow properties and dependence on Riemannian metrics
are discussed in \S 5. The proof of Theorem A (Proposition~\ref{5.2}
for (1), Proposition~\ref{corb} for (2) and Proposition~\ref{6.5} for
(3)) is occupied in \S 6 and \S 7. In \S 8, we study the relative Seiberg-Witten
invariant and complete
the proof of Theorem B as Theorem~\ref{invariant}.

\section{Seiberg-Witten equation on 3-manifolds}

It is well-known that every closed oriented 3-manifold is spin. 
The group $Spin(3) \cong SU(2) \cong Sp_1$
is the universal covering of $SO(3) = Spin(3)/\{\pm I \}$.
Pick a Riemannian metric $g$ on $Y$. The 
metric $g$ defines the principal $SO(3)$-bundle $P_{SO}(Y)$ of oriented
orthonormal frames on $Y$. A spin structure is a lift of
$P_{SO}(Y)$ to a principal $Spin(3)$-bundle $P_{Spin}(Y)$ over $Y$. The set
of equivalence classes of such lifts has, in a natural way, the structure of a
principal $H^1(Y, Z_2)$-bundle over a point. So there is a unique
spin-structure on the integral homology 3-sphere $Y$.
\medskip

There is a natural adjoint representation
\[Ad: Spin(3) \x Sp_1 \to Sp_1; \ \ \ (q, \a ) \mapsto q \a q^{-1}, \]
and associated rank-2 complex vector bundle (spinor bundle)
\[W = P_{Spin(3)}(Y) \x_{Ad} C^2 .\] 
Let $L = det W$ be the determinant line bundle. 
For the ordinary Spin-structure,
one has a Clifford multiplication
\[c: T^*Y \ox W \to W\]
\[ c([p,\a]) \ox [p,v] \to [p, \overline{\a}v].\]
So $c$ induces a map $T^*Y \to Hom(W, W)$. The spinor pairing 
$\t : W \ox \overline{W} \to T^*Y$ is given by
\[[p, v_1 \ox v_2] \to \t(\frac{1}{4} Im(v_1 i v_2)),\]
where $\t$ is an orientation preserving isomorphism $P_{Spin(3)}(Y) \x
Sp_1 \to T^*Y$. A connection $a$ on $L$ together with the Levi-Civita connection
on the tangent bundle of $Y$ form a covariant derivative on $W$.
This maps sections of $W$ into sections of $W \ox T^*Y$. Followed
by the Clifford multiplication, one has a Dirac operator
\[{\bd }^g_a : \G (W) \stackrel{\nabla^g_a}{\to } 
\G (W \ox T^*Y) \stackrel{c}{\to } \G (W). \]
The determinant line bundle $L$ is trivial for the spin structure, so we may
choose $\theta$ to be the trivial connection and ${\bd }^g_{\theta }:
\G (W) \to \G (W)$ is the usual Dirac operator. Note that all bundles over the 
integral homology 3-sphere $Y$ are {\bf trivial}.

There is a unique spin-structure on $Y \x \R$ associated to the unique
spin-structure on $Y$ with the product metric on $Y \x \R$.
The two spinor bundles $W^{\pm }$ on $Y \x \R$ can be identified 
by using a Clifford
multiplication by $dt$, where $t$ is denoted for the variable on $\R$.
Both $W^+$ and $W^-$ are obtained by the pull-back of the
$U(2)$-bundle $W \to Y$ from the projection map $Y \x {\R} \to Y$.
Thus we have the identification of the map $\s :
{\Lam }^2T^*(Y \x {\R}) \to Hom(W^+, W^-)$ and the map
$\t^{-1} : T^*Y \to Hom(W, W)$ through the above
identifications. $\s (\e ) = \t^{-1}(*_g \e)$. In other words from
the identification ${\Lam }^2T^*(Y \x {\R}) = {\Lam }^2T^*Y \oplus 
{\Lam }^1T^*Y$ and using the Hermitian pairing on $W^{\pm }$, there is an
induced pairing
\[\t : \overline{W} \x W \to {\Lam }^1T^*Y .\]
In fact for every $\g: T^*Y \to Hom(W, W)$ (a spin structure),
that is a way to determine a spin structure on $Y \x \R$ by 
\[\s :T^*(Y \x {\R}) \to Hom(W \oplus W, W \oplus W); \ \ \
\s (v, r) = \left( \begin{array}{cc}
0 & \g (v) + r 1 \\ \g (v)  - r 1 & 0 \end{array} \right) .\]
The determinant line bundle $L_{(4)} = det W^{\pm }|_{Y \x \R}$ (a
trivial line bundle) carries $U(1)$-connections $A = a + \p dt$. So the Dirac
operator $D^g_A$ for the product metric
$g + dt^2$ over $Y \x \R$ is given by
\[ D^g_A = \left( \begin{array}{cc}
0& - \nabla_t + {\bd }^g_a \\
\nabla_t + {\bd }^g_a & 0 \end{array} \right) , \]
where ${\bd }^g_a$ is a twisted self-adjoint Dirac operator on 
$\G (W) \to \G (W)$, and $\nabla_t = \frac{\bd }{\bd t} + \p$ is a twisted
skew adjoint Dirac operator
over $\R$.

The curvature 2-form of $A = a + \p dt$ can be calculated as 
$F_A = F_a + (\frac{\bd a}{\bd t} - d_a \p) dt$.
Using the identification of ${\O }^2(Y \x {\R}) \cong {\O }^2(Y) \oplus
{\O }^1(Y)$, we can write $F^+_A$ as
$*_g F_a + (\frac{\bd a}{\bd t} - d_a \p) \in {\O }^1(Y)$
as the self-dual component of the curvature $F_A$. Now the 
Seiberg-Witten monopole equation on 4-manifolds
reduces to a Seiberg-Witten monopole equation on 3-manifolds as
\begin{equation} \label{SW}
\left\{ \begin{array}{ll}
(\nabla_t + {\bd }^g_a) \psi & = 0 \\
*_g F_a + (\frac{\bd a}{\bd t} - d_a \p)& = i \t(\psi, \psi) 
\end{array}  \right. 
\end{equation}
for $\psi \in \G (W)$. It is equivalent to the flow equation of $(a, \p, 
\psi )$:
\begin{equation} \label{SW0}
\left\{ \begin{array}{ll}
\frac{\bd \psi }{\bd t} & = - {\bd }^g_a \psi - \p . \psi \\
\frac{\bd a}{\bd t} & = - *_g F_a + d_a \p + i \t(\psi, \psi) .
\end{array}  \right. 
\end{equation}
The equation (\ref{SW}) is invariant under the gauge transformation
$u \in Map (Y, U(1))$, where the gauge group action on
$(a + \p dt, \psi )$ is given by
\begin{equation} \label{star}
u \cdot (a + \p dt, \psi ) = (u^*a + 
(\p - u^{-1}\frac{du }{dt})dt, \psi u^{-1}).
\end{equation}
There is a temporal gauge to obtain a simpler equation. The temporal gauge $u$
is the element which $u \cdot (a + \p dt) = u^*a$, i.e.,
$\p - u^{-1}\frac{du }{dt} = 0$. Then the equation (\ref{SW0}) can be reduced
to the following form.
\begin{equation} \label{rSW}
\left\{ \begin{array}{ll}
\frac{\bd \psi }{\bd t} & = - {\bd }^g_a \psi \\
\frac{\bd a}{\bd t} & = - *_g F_a +i \t(\psi, \psi) .
\end{array}  \right. 
\end{equation}

\section{Configuration spaces on $Y$}

Fix a trivialization $L = Y \x U(1)$, one can identify the
space of $U(1)$-connections of Sobolev $L_k^p$-norm with the space
${\CA}_k^p = L_k^p({\O }^1(Y, i \R))$ of 1-forms on $Y$ such that the 
zero element in ${\O }^1(Y, i \R)$ corresponds to the trivial
connection $\theta $ on $L$. The gauge group of $L$ can be identified
with ${\CG }^p_k(Y) = L^p_{k+1}(Map(Y, U(1)))$ acting on
${\CA}_k^p \x L^p_k(\G (W))$ by (\ref{star}). We need to assume that
$k+1 > 3/p$ so that ${\CG}_Y = {\CG }^p_k(Y)$ is a Lie group. We may take
$k=1, p=2$. 

Let ${\CC}_Y$ be the configuration space
\[{\CC}_Y = L^2_k(\{{\O }^1 \oplus {\O }^0\}(Y, i {\R}) \oplus \G (W)) .\]
The quotient space is ${\CB}_Y = {\CC}_Y/{\CG }_Y.$ Denote
${\CC}_Y^* = \{(a, \p, \psi ) \in {\CC}_Y| \psi \neq 0 \}$.
For $(a, \p, \psi ) \in {\CC}_Y^*$, the isotropy group ${\G }_{(a, \p, \psi )} 
= \{id \}$. For $(a, \p, \psi ) \in {\CC}_Y \setminus {\CC}_Y^*$, the isotropy
group ${\G }_{(a, \p,0)} = U(1)$, these elements are called reducibles. 
For example, $\T = (\theta , 0, 0)$ is reducible by all constant
maps from $Y$ to $U(1)$. Note that ${\CG }_Y$ acts
freely on ${\CC}_Y^*$, so
${\CB}_Y^* = {\CC}_Y^*/{\CG }_Y$ forms an open and dense set in
${\CC}_Y/{\CG }_Y$.
\begin{pro} \label{hil}
${\CB}_Y^*$ is a Hilbert manifold. For $(a_0, \p_0, \psi_0) \in {\CC}_Y^*$,
the tangent space of ${\CB}_Y^*$ can be identified with
\[T_{[(a_0, \p_0, \psi_0)]}{\CB}_Y^* = \{
(a, \p, \psi ) \in 
L^2_k(\{{\O }^1 \oplus {\O }^0\}(Y, i {\R}) \oplus \G (W))| \]
\[ \|(a, \p, \psi )\|_{L^2_{k-1}} < \ve, \ \ \ 
d^*_{a_0}\psi + Im (\psi_0 , \psi ) = 0 \}. \]
\end{pro}
\noindent{\bf Proof:}
This follows from the construction of slice in \cite{dk, fu}.
It will be clear from context to identify $(a_0, \p_0, \psi_0)$ with its
gauge equivalence class in our notation. The gauge orbit of 
$(a_0, \p_0, \psi_0) \in {\CC}_Y^*$ is given by ${\CG }_Y \to {\CC}_Y^*$:
\[g=e^{iu} \to (a_0 - g^{-1} dg , \p_0 , \psi_0 g^{-1}) .\]
The linearization of this map at $Id = e^0$ is
\[\d_0 : T_{id} {\CG }_Y = {\O }^0(Y, i{\R}) \to 
\{{\O }^1 \oplus {\O }^0\}(Y, i {\R}) \oplus \G (W) \]
\[ u \mapsto (-du, 0, - \psi_0 u) .\]
So the adjoint operator $\d_0^*$ of $\d_0$ is given by
\[ \d_0^* \psi = d^*_{a_0} \psi + Im (\psi_0 . \psi ) .\]
A neighborhood of $[(a_0, \p_0, \psi_0)] \in
{\CB}_Y^*$ can be described as a quotient of 
$T_{[(a_0, \p_0, \psi_0)], \ve}{\CB}_Y^* /{\G }_{(a_0, \p_0, \psi_0)}$
for sufficiently
small $\ve $. Every nearby orbit meets the slice
$(a_0, \p_0, \psi_0) + T_{[(a_0, \p_0, \psi_0)], \ve}{\CB}_Y^*$. This is amount
to solving the gauge fixing condition relative to
$(a_0, \p_0, \psi_0)$, i.e., there exists a unique
$u \in {\O }^0(Y, i {\R})$ such that
$e^{iu} \cdot (a_0 + a, \p_0 + \p, \psi_0 + \psi ) \in 
T_{[(a_0, \p_0, \psi_0)], \ve}{\CB}_Y^*$
for $\psi_0 \neq 0$. Hence it follows from applying the implicit function 
theorem. \qed

There is an associated bundle ${\CC}_Y^* \x_{{\CG }_Y}({\O }^1(Y, i {\R})
\oplus \G (W))$ over ${\CC}_Y^*$ because of the free action of
${\CG }_Y$ on ${\CC}_Y^*$. We define a section $f :
{\CC}_Y^* \to {\CC}_Y^* \x_{{\CG }_Y}({\O }^1(Y, i {\R})
\oplus \G (W))$ by
\[f(a, \p , \psi ) = [(a, \p , \psi ),
*_gF_a -d_a \p - i \t (\psi, \psi ), {\bd }^g_a \psi + \p . \psi ]. \]
Note that $f$ is ${\CG }_Y$-equivariant, 
$f(g \cdot (a, \p , \psi )) = g \cdot 
f(a, \p , \psi )$. Hence it descends to ${\CB}_Y^*$,
\[f : {\CB}_Y^* \to {\CC}_Y^* \x_{{\CG }_Y}({\O }^1(Y, i {\R}) \oplus \G (W)).\]
Now $f(a, \p , \psi ) \in T_{[(a, \p , \psi )], \ve}L^2_{k-1}
{\CB}_Y^* = {\CL }_{[(a, \p , \psi )]}$. So $f$ can be thought of as
a vector field on the Hilbert manifold ${\CB}_Y^*$. Over ${\CB}_Y^*$, $f$ is 
a section of the bundle ${\CL }$ with fiber ${\CL }_{[(a, \p , \psi )]}$.
\begin{df} \label{flat}
The zero set of $f$ in ${\CB}_Y^*$ is the moduli space of solutions of
the 3-dimensional Seiberg-Witten equation
\[f^{-1}(0) = {\CR}_{SW}^*(Y, g) = \{
[(a, \p , \psi )] \in {\CC}_Y^* \ \mbox{satisfies (\ref{SWF})}\}/{\CG }_Y.\]
\begin{equation} \label{SWF}
\left\{ \begin{array}{l}
{\bd }^g_a \psi + \p . \psi = 0 \\
*_g F_a - d_a \p - i \t (\psi, \psi ) = 0
\end{array} . \right.
\end{equation}
\end{df}
We will show that ${\CR}_{SW}^*(Y, g)$ is a zero-dimensional smooth
manifold and its algebraic number is the Euler characteristic of a monopole
homology defined in \S 6 (see also \cite{au} for instance). 

The linearization of $f$ can be computed as the following.
\begin{eqnarray*}
f(a_0 + sa, \p_0 + s \p , \psi_0 + s \psi ) & = & 
(*_g F_{a_0 + sa} - d_{a_0 + sa}(\p_0 + s \p) - i \t (\psi_0 + s \psi ,\\
& &
\psi_0 + s \psi ), {\bd }^g_{a_0 + sa}(\psi_0 + s \psi ) +
(\p_0 + s \p ). (\psi_0 + s \psi ) \\
& = & f(a_0, \p_0, \psi_0)  +s \d_1(a_0, \p_0, \psi_0)((a, \p , \psi )) +
o(s^2) .
\end{eqnarray*}
So the linearized operator $Df(a_0, \p_0, \psi_0) = 
\d_1(a_0, \p_0, \psi_0) : T_{[(a_0, \p_0, \psi_0)]} {\CB }^*_Y
\to {\CL }_{[(a_0, \p_0, \psi_0)]}$
is given by 
\[\d_1(a_0, \p_0, \psi_0) : \{{\O }^1 \oplus {\O}^0\}(Y, i{\R })\oplus \G (W)
\to {\O }^1(Y, i{\R })\oplus \G (W) ,\]
\[((a, \p , \psi ) \longmapsto \left( \begin{array}{ccc}
*_g d_{a_0} & - d_{a_0} & - i Im(\psi_0,, \cdot) \\
c(\cdot \psi_0) & c \cdot \psi_0 & {\bd }^g_{a_0} + \p_0 \cdot 
\end{array} \right) 
\left(\begin{array}{c}
a \\ \p \\ \psi 
\end{array} \right) .\]
It forms a natural 3-dimensional monopole complex, since 
$\ker \d_0^*$ is the gauge fixing slice. So
\begin{equation} \label{mc}
MC_{\bullet}:  0 \to {\O}^0(Y, i{\R }) \stackrel{\d_0}{\to } 
\{{\O }^1 \oplus {\O}^0\}(Y, i{\R })\oplus \G (W) 
\stackrel{\d_1}{\to } {\O }^1(Y, i{\R })\oplus \G (W) \to 0,\end{equation}
is a short exact sequence. The operator
\[\d_0^* \oplus \d_1 (a_0, \p_0, \psi_0): 
\{{\O }^1 \oplus {\O}^0\}(Y, i{\R })\oplus \G (W)
\to \{{\O }^1 \oplus {\O}^0\}(Y, i{\R })\oplus \G (W) \]
\begin{equation} \label{32*}
(a, \p , \psi ) \longmapsto \left( \begin{array}{ccc}
*_g d_{a_0} & - d_{a_0} & - i Im(\psi_0, \cdot) \\
-d_{a_0}^* & 0 & Im (\psi_0, \cdot ) \\
c(\cdot \psi_0) & c \cdot \psi_0 & {\bd }^g_{a_0} + \p_0 \cdot 
\end{array} \right) 
\left(\begin{array}{c}
a \\ \p \\ \psi 
\end{array} \right) ,\end{equation}
is a first-order operator with symbol
$\s (\d_0^* \oplus \d_1) = \s (\d)$, where
\[ \d = \left( \begin{array}{ccc}
*_g d_{a_0} & - d_{a_0} & 0 \\
- d_{a_0}^* & 0 & 0 \\
0 & 0 & {\bd }^g_{a_0}
\end{array} \right) \]
is a first-order self-adjoint Dirac operator. Hence
\begin{eqnarray}
Ind (\d_0^* \oplus \d_1) & = & Ind (\d) \nonumber \\
 & = & Ind \left( \begin{array}{cc}
*_g d_{a_0} & - d_{a_0} \nonumber  \\
- d_{a_0}^* & 0 \end{array} \right) + Ind {\bd }^g_{a_0} \nonumber \\
 & = & 0 .
\end{eqnarray}
Since  the operator $\left( \begin{array}{cc}
*_g d_{a_0} & - d_{a_0} \\ - d_{a_0}^* & 0 \end{array} \right)$ is self-adjoint
and every Dirac operator has index zero over odd (3-)dimensional manifolds, thus
we have the zero index for the operator $\d_0^* \oplus \d_1$. 
Generically, the moduli space 
${\CR }_{SW}(Y, g)$ is zero-dimensional. 

Define $H^0(MC_{\bullet}) = \ker \d_0$, $H^1(MC_{\bullet}) = \ker \d_1/im \d_0$,
$H^2(MC_{\bullet}) = coker \d_1$. The first cohomology
$H^1(MC_{\bullet})$ is isomorphic for 
every $(a_0, \p_0, \psi_0) \in {\CB }^*_Y$, 
so that $(a_0, \p_0, \psi_0) \in {\CB }^*_Y$ is a nondegenerate zero 
of $f$ if and only if $\ker (\d_0^* \oplus \d_1) = H^1(MC_{\bullet}) = 0$.
For $\T = (\theta, 0, 0)$ and a generic metric $g$ without harmonic spinors of
$\bd_{\theta}^g$, we have that $\T$ is always isolated and nondegenerate (in the
Bott sense) zero of $f$ on the integral homology 3-sphere $Y$.

\section{Admissible Perturbation and Transversality}

In this section, we prove that there are enough perturbations
to make the zero set of $f$ transverse.
There is a 1-form perturbation reduced from 4-dimensional Seiberg-Witten
equation as in \cite{don, km, tau}. In our 3-dimensional
case, {\em the harmonic spinor may vary or jump as metrics on $Y$ vary}. 
In order to obtain any topological information, one needs to 
extend the perturbation-data and understand the harmonic spinors
accordingly.
The method we used here is essentially the one used in \cite{fu, ll, li, pt}.

Let ${\CP}_Y = \Si_Y \x {\O}^1 (Y, i {\R})$ be the space of perturbation
data, where $\Si_Y$ is the space of Riemannian metrics on $Y$. Consider
the union 
$\cup_{(g,, \a) \in {\CP}_Y}{\CR}^*_{SW}(Y;g, \a)$ of the moduli spaces
of  3-dimensional Seiberg-Witten solutions over all metrics and 1-forms. If
the union is a (Banach) Hilbert manifold, then its projection to the space
${\CP}_Y$ is a Fredholm map. So there exists a Baire first  category in
${\CP}_Y$ such that ${\CR}^*_{SW}(Y;g, \a)$ is a manifold by the
Sard-Smale theorem.

Let $f_{\e}$ be the parametrized smooth section of the bundle
${\CL} \to {\CB}^*_Y \x {\CP}_Y$ with 
$\e = (g, \a ) \in {\CP}_Y$. The map $f_{\e}$
is given by
\[ f_{\e} : {\CB}^*_Y \to {\O}^1(Y, i {\R}) \oplus \G (W) \]
\[ (a, \p, \psi ) \mapsto (*_gF_a - d_a \p - i \t (\psi , \psi) + \a ,
{\bd }^{\nabla_0 + \a}_{a} \psi + \p . \psi ) , \]
where $\nabla_0$ is the Levi-Civita connection for the metric $g$.
Let $f_{1 \e}(a, \p, \psi ) = {\bd }^{\nabla_0 + \a}_a \psi + \p . \psi$
be the second component of the map $f_{\e }$ on $\G (W)$, and
$f_{0 \e}(a, \p, \psi )$ be the first component of $f_{\e}$.

\begin{lm} \label{sub}
$f_{1 \e}$ is a submersion ($Df_{1 \e}$ is surjective).
\end{lm}
\noindent{\bf Proof:} The differential $Df_{1 \e}$ is given by the formula
\[Df_{1 \e}(a, \p, \psi ; o , \a)(\ve a , \ve \p , \ve \psi , 0, \ve \a)
= {\bd }^{\nabla_0 + \a }_a (\ve \psi ) +
(\ve \a + \ve a + \ve \p ). \psi + \p . \ve \psi , \]
where we vary along the subspace 
$\{{\O}^1 \oplus {\O}^0(Y, i {\R}) \oplus \G (W)\} \x 
\{\{0\} \x {\O}^1(Y, i {\R})\}$ of $T^*_{[a, \p, \psi]}{\CB }^*_Y \x
{\CP}_Y$. We want to show that
$Df_{1 \e}$ is surjective. Suppose the contrary. Then there exists a spinor
$\c \in \G (W)$ such that it is perpendicular to $Im Df_{1 \e}$. 
\begin{equation} \label{*a}
\la {\bd }^{\nabla_0 + \a }_a (\ve \psi ) , \c \ra =0, 
\end{equation}
for all $\ve \psi$. I.e., $\c \in \ker ({\bd }^{\nabla_0 + \a }_a)^*$.
By the elliptic regularity of (\ref{*a}), a solution $\c$ is smooth. Choose
a point $y \in Y$ such that $\c (y) \neq 0$. By the uniqueness of continuation 
of the solution of the elliptic equation \cite{an}, 
${\bd }^{\nabla_0 + \a }_a \cdot ({\bd }^{\nabla_0 + \a }_a)^* \c = 0$,
there is a neighborhood $U_y$ of $y$ such that $\c (y) \neq 0$ for $y \in U_y$.
Thus we can find a 1-form $\ve \a + \ve a \in {\O }^1(Y, i{\R})$ such that
$(\ve \a + \ve a). \psi = \lam \c$ with $\lam \neq 0$ in $U_y$, and
$\ve \a + \ve a$ has compact support. So we obtain
\begin{eqnarray*}
0 & = & \la {\bd }^{\nabla_0 + \a + \ve \a}_{a+ \ve a} (\ve \psi ) , \c \ra \\
 & = & \la {\bd }^{\nabla_0 + \a }_a (\ve \psi ) , \c \ra +
\la (\ve \a + \ve a). \ve \psi , \c \ra \\
 & = & \la \lam \c , \c \ra 
 = \lam \la \c , \c \ra .
\end{eqnarray*}
Therefore $\c = 0$ in $U_y$, so
$\c \equiv 0$ by a result in \cite{an}. \qed

By the Hodge decomposition of ${\O }^1(Y, i{\R}) = Im d \oplus Im d^*$ for
$Y$, we have that $\d_1$ is surjective. Thus
$f_{0\e }(\a , \p, \psi ) = *_g F_a - d_a \p - i \t (\psi, \psi ) + \a $
is also a submersion onto ${\O }^1(Y, i{\R})$. 

\begin{cor} 
The spaces $f_{0\e }^{-1}(0)$ and 
$f_{1 \e }^{-1}(0)$ are Banach manifolds.
\end{cor} \qed

Now at point $(a_0, \p_0, \psi_0; g_0, \a) \in {\CC}_Y \x {\CP}_Y$, 
the parametrized smooth section 
\[f(a_0, \p_0, \psi_0; g_0, \a) = f_{(g_0, \a)}(a_0, \p_0, \psi_0) =
f_{\e }(a_0, \p_0, \psi_0 ) \]
is submersion.

\begin{pro} 
The differential $Df$ is onto at all points of the 
moduli space $f^{-1}(0) \subset {\CC}_Y^* \x {\CP}_Y$.
\end{pro}
\noindent{\bf Proof: }
The differential $Df$ at $(a_0, \p_0, \psi_0; g_0, \a) \in {\CC}_Y \x {\CP}_Y$
is of the form $(Df_0, Df_1)$
\begin{eqnarray*}
Df_0 & = & *_{g_0}d_{a_0}a + (g)_* F_{a_0} - d_{a_0} \p - i Im (\psi_0, \psi )
- a . \p_0 + \a \\
Df_1 & = & {\bd }^{\nabla_0 + \a_0}_{a_0} \psi + (\a + a). \psi_0
+ (\p . \psi_0 + \p_0 . \psi ) + r(g) )
\end{eqnarray*}
where $(g)_*$ is the variation of the Hodge star operator
$(g)_* = \frac{d}{ds}|_{s =0} *_{g_0 + sg}$, $r(g)$ is a zero order
operator applied to the variation $g_0 + sg + o(s^2)$ of metric, 
$a. \p_0$ is the Clifford multiplication of 1-form $a$ on the section
$\p_0 \in \G(W)$. The surjective of 
$Df_0$ follows from Theorem 3.1 of \cite{fu}, and
the surjective of $Df_1$ follows from Proposition I.3.5 of \cite{pt}
(see also \cite{au, m, mt, moy}). 
\qed

We consider the map $f_*: {\CC}_Y^* \x {\CP}_Y \to 
{\O}^1(Y, i{\R}) \oplus \G (W)$.
\begin{cor} 
The space $f_*^{-1}(0)$ is a Banach manifold.
\end{cor}
\noindent{\bf Proof:} Take $f_*$ as a section of ${\CB}^*_Y \x {\CP}_Y$ to 
$({\CC}_Y^* \x_{{\CG}_Y} ({\O}^1(Y, i{\R }) \oplus \G (W)) \x {\CP }_Y$.
So $f_*^{-1}(0)|_{{\CB}^*_Y} = f_*^{-1}(0)/{\CG}_Y$ is a Banach manifold.
\[\begin{array}{ccc}
{\CC}_Y^* \x {\CP }_Y & \stackrel{f}{\to } & {\O}^1(Y, i{\R }) \oplus \G (W) \\
 \downarrow \pi_2 & & \\
{\CP }_Y & &  
\end{array} \]
The projection map $\pi_2$ is a smooth Fredholm map
of index zero. It follows exactly from the same argument in \cite{dk, fu}. 
\qed
\begin{cor} \label{admi}
The inverse image $\pi_2^{-1}((g, \a))$ of a generic
parameter $(g, \a) \in {\CP }_Y$, 
the moduli space ${\CR }_{SW}(Y, (g, \a ))$ of
the 3-dimensional monopole solutions is a zero
dimensional manifold.
\end{cor}

A perturbation $\e = (g, \a)$ satisfying Corollary~\ref{admi} is called
{\bf admissible.}
In general, the class of reducible elements in ${\CC}_Y \setminus {\CC}_Y^*$ 
forms
a singular strata in the quotient space ${\CB}_Y$. If it is a solution
of 3-dimensional Seiberg-Witten equation, it is also singular to the space
of ${\CR}_{SW}(Y, g)$. The reducible solutions
of the 3-dimensional Seiberg-Witten equation satisfy
\begin{eqnarray} \label{reds}
{\bd }^{\nabla_0 + \a}_{a} \psi + \p_0 . \psi & = & 0 \nonumber \\
 - *F_a + d_a \p & = & 0, 
\end{eqnarray}
for $\psi = 0$. Applying the temporal gauge $g \cdot (a, \p ) = (g^*a, 0)$, 
we get that $g^*a$ is a flat connection on $Y \x U(1)$ over $Y$. For integral
homology 3-sphere, there is a unique $U(1)$ reducible connection, 
namely the trivial one. So the reducible solution is $(\theta , 0)$. 
There is a unique $U(1)$-reducible solution of (\ref{reds}), denoted
by $\Theta = (\theta, 0)$. 

Note that
$\ker \d_1 = \ker {\bd }^{g}_a$ for an integral homology 3-sphere. For
a generic metric $g$, $\ker {\bd }^{g}_a = 0$. But $\ker {\bd }^{g_t}_a$
may have a nontrivial kernel as the Riemannian metrics vary in an
one-parameter family (see \cite{hi}).  The harmonic spinor, even the dimension
of the harmonic spinor, depends on the metric used in defining the Dirac
operator. Hence
the harmonic-spinor jump creates and/or destroys irreducible solutions
of the 3-dimensional Seiberg-Witten equation.
This is the main problem to understand the new phenomenon that the
``Seiberg-Witten-Floer theory'' is not entirely
metric-independent (see \cite{don}).
In the next section, we study such a dependence of Riemannian metrics.

\begin{pro} \label{path}
${\CR }_{SW}^*(Y, (g, \a)) =  {\CR }_{SW}(Y, (g, \a)) \setminus \{ \Theta\}$
is a zero-dimensional smooth manifold for a first category near $(g, \a)$ in
${\CP}_Y$.
\end{pro}
\noindent{\bf Proof:} The results follows from the construction above,
Proposition 2c.1 of \cite{fl} and the Sard-Smale theorem. 
\qed

Define the weighted Sobolev space $L_{k,\d}^p$ on sections $\xi $
of a bundle over $Y\x {\R}$ to be the space of $\xi$ for
which $e_{\d}\cdot \xi $ is in $L_k^p$,
where $e_{\d}(y,t) = e^{\d |t|}$ for
$|t| \geq 1$. For any $\d \geq 0$ and any Seiberg-Witten monopole
solution $(A, \Phi)$ on $Y \x {\R}$, the linearized operator
\[D_{A, \Phi}: L_{k+1,\delta }^p(\G (W_{(4)}^+) \oplus
\O^1(Y \x {\R})) \to L_{k,\delta }^p(\G (W_{(4)}^-) \oplus
({\Omega }^0 \oplus {\Omega }^2_+)(Y \x {\R}))\]
is Fredholm (see \cite{don, fl, km, tau, wi}).
We call $(A, \Phi)$ {\it regular} if $\mbox{Coker} D_{A, \Phi} = 0$
and we call ${\CM}_{Y \x {\R}}$ (the moduli space of
perturbed Seiberg-Witten solutions with finite energy) {\it regular}
if it contains orbits of regular $(A, \Phi)$'s.

\begin{pro} \label{bfc}
The finite energy condition forces elements of ${\CM}_{Y \x {\R}}$
to converge to zeros of $f^{-1}_{\e}(0)$ on the ends of
$Y\x {\R}$. The set of all perturbations $\e \in {\CP}_Y$ of which
${\CM}_{Y \x {\R}}$ is regular is of Baire's first category.
\end{pro}
\noindent{\bf Proof:} 
The proof follows exactly from the same method in \cite{fl} Proposition
2c.2 with Chern-Simons Seiberg-Witten functional
as defined in \cite{km} \S 4 and \cite{au, m, mt}. \qed

\section{Spectral flow and Dependence on Riemannian metrics}

In this section, we use the unique $U(1)$-reducible
solution $\Theta$ to capture the metric-dependent relation via the spectral
flow. In \cite{ll} joined with Lee, the author
used the Walker correction-term around 
$U(1)$-reducibles to obtain homotopy classes of admissible perturbations
(realized by a family of Lagrangians), and to show the invariance among
the same homotopy class of the Lagrangian perturbations. Those
Walker correction-term can be interpreted as the spectral flow 
in \cite{clm, ll}.

\begin{pro} \label{sfi}
For an admissible perturbation $\e = (g, \a) \in {\CP}_Y$ and a nondegenerate
zero $(a, \phi, \psi) \in {\CR}_{SW}(Y, \e) = f_{\e}^{-1}(0)$,
we can associate an integer $\mu_{\e}(a, \phi, \psi) \in {\Z}$ such that
for $(A, \Phi) \in {\CB}_{Y \x {\R}}((a, \phi, \psi),
(a^{'}, \phi^{'}, \psi^{'}))$
\begin{eqnarray*}
\mu_{\e}(e^{iu} \cdot (a, \phi, \psi)) &= &\mu_{\e}(a, \phi, \psi), \\
\mbox{Index}D_{A, \Phi} & = & \mu_{\e}(a, \phi, \psi) - 
\mu_{\e}(a^{'}, \phi^{'}, \psi^{'}) - \mbox{dim} 
\G_{(a^{'}, \phi^{'}, \psi^{'})},
\end{eqnarray*}
where $\G_{(a^{'}, \phi^{'}, \psi^{'})}$ is the isotropy subgroup of
$(a^{'}, \phi^{'}, \psi^{'})$.
\end{pro}
\noindent{\bf Proof:}
Let $\pi_1: Y \x [0, 1] \to Y$ be the projection on the first factor.
Let $L_{(4)} \x W_{(4)}$ be the pullback 
$\pi_1^*(det W^{\pm}) \x \pi_1^*W^{\pm}$ such that 
$(A, \Phi) \in {\CA}_{L_{(4)}} \x W_{(4)}$ satisfies
$(A, \Phi)|_{t \leq 0} = (a, \phi, \psi)$ and
$(A, \Phi)|_{t \geq 1} = (a^{'}, \phi^{'}, \psi^{'})$. We have
$D_{A, \Phi} = \frac{\bd}{\bd t} + \d_t$ with $\d_t = \d_{A(t), \Phi(t)}$
in (\ref{32*}). Then the Fredholm index of $D_{A, \Phi}$ is given by the
spectral flow of $\d_t$ (see \cite{aps, clm, fl}). The second equality
follows from the same proof of Proposition 2b. 2 in \cite{fl}. The first
equality follows from
\begin{eqnarray*}
SF(e^{iu}\cdot (a, \phi, \psi), (a, \phi, \psi))& = &
\mbox{Ind}D_{A, \Phi}((a, \phi, \psi), (a, \phi, \psi))_{Y\x S^1} \\
& = & \frac{1}{4}(c_1(L_{(4)})^2 - (2 \c + 3 \s))(Y\x S^1) = 0,
\end{eqnarray*}
where $\c$ and $\s$ are the Euler number and signature of $Y \x S^1$, and
$c_1(L_{(4)})^2(Y\x S^1) = 0$ for the integral homology 3-sphere $Y$.
\qed

Note that the relative index is gauge-invariant, but depending on the
perturbation $\e \in {\CP}_Y$ by Proposition~\ref{sfi}.
The absolute index may not be well-defined since $\mu_{\e}(\Theta)$
depends upon $\e \in {\CP}_Y$. In the instanton case, we fix the trivialization
of a principal bundle
and a fixed tangent vector to the trivial connection to determine
$\mu (\theta) =0$ for the trivial connection $\theta$. It turns out
that such a fixation is independent of metrics and other perturbation data
in the instanton Floer theory.
But this is no longer true for the monopole case.

\begin{pro} \label{1pf}
{\bf (Definition)}
Two admissible perturbations $\e_0$ and $\e_1$ in ${\CP}_Y$ are (called) 
homotopic
to each other through a 1-parameter family $\e_t (0 \leq t \leq 1)$ in
${\CP}_Y$ if and only if $\mu_{\e_0}(\T) = \mu_{\e_1}(\T)$.
\end{pro}
\noindent{\bf Proof:}
For two admissible perturbations $\e_0$ and $\e_1$ in \S 4, we can
connect them into a 1-parameter family $\e_t$ such that there are
at most finitely many
$t\in (0, 1)$ with $\e_t$ corresponding harmonic-spinor jumps. Denote those
$0< t_0 < t_1  \cdots < t_n < 1$ and $\lam_1, \lam_2, \cdots, \lam_n,
\lam_{n+1} = 0$ so that $\lam_i$ is not the eigenvalues of $\d_t = 
\d_t (\theta, 0)$ for $t_{i-1} \leq t \leq t_i$, where
$t_{-1 } =0$ and $t_{n+1} =1$.
Define $n_i = \mbox{dim}(\d_{t_i} - \lam Id)$ with $\lam \in 
[\lam_{i+1}, \lam_i]$ and $n_i = - \mbox{dim}(\d_{t_i} - \lam Id)$ with 
$\lam \in [\lam_i, \lam_{i+1}]$. From the operator
$D_{\e_t}(\Theta) = \frac{\bd}{\bd t} + \d_t(\Theta)$ and the well-known facts
in \cite{aps, clm, fl}, we have
\[\mbox{Ind} D_{\e_t}(\Theta) = \sum_{i=0}^n n_i.\]
This shows that $\mbox{Ind} D_{\e_t}(\Theta)$ is independent of
the construction $\e_t$ and that is continuous in $\e_t$.
On the other hand, \[\mbox{Ind} D_{\e_t}(\Theta) =
\mu_{\e_0}(\Theta) - \mu_{\e_1}(\Theta).\]
Thus the obstruction to connect two generic perturbations is the
spectral flow along the metric path in $\Si_Y$. The 
Riemannian-metric space $\Si_Y$ is path-connected.
So $\mbox{Ind} D_{\e_t}(\Theta) =0$ provides that
$\e_0$ and $\e_1$ are in the same (homotopy) class of 
with respect to the spectral flow.
\qed

Thus the dependence of metrics also enters into the definition of 
relative indices for $(a, \phi, \psi) \in {\CR}_{SW}^*(Y, \e)$. Now we follow
the instanton case to fix the relative index
\[\mu_{\e}(a, \phi, \psi) = 
\mbox{Ind}D_{\e}(\T, (a, \phi, \psi)) \in {\Z}, \]
which depends on the value $\mu_{\e}( \Theta)$. Any changes of
$\mu_{\e}( \Theta)$ shift $\mu_{\e}(a, \phi, \psi)$ by an integer,
and $\mu_{\e}(\T)$ is understood with respect to some reference
perturbation $\e_0 \in {\CP}_Y$.

\begin{lm} \label{cod}
For an admissible perturbation $\e \in {\CP}_Y$, the Seiberg-Witten moduli space
${\CR}_{SW}(Y, \e) = f_{\e}^{-1}(0)$ is a
compact 0-dimensional oriented manifold.
\end{lm}
\noindent{\bf Proof:}
The compactness can be proved by the 3-dimensional Weitzenb\"{o}ck
formula and Mosers' weak maximal principle as in the 4-dimensional case
\cite{km, wi}. By the construction in the proof of Proposition~\ref{sfi}, 
we can show that ${\CR}_{SW}(Y, \e) = f_{\e}^{-1}(0)$ is a closed
subset of the compact moduli space 
${\CM}_{Y\x S^1}(g + d\theta, \pi_1^* \e)$, where $Y \x S^1$
carries the product metric $g + d\theta$. That
${\CR}_{SW}(Y, \e)$ is compact follows by Lemma 2 of \cite{km}.
By Proposition~\ref{path}, ${\CR}_{SW}(Y, \e)$ is a
0-dimensional manifold. The orientation at each point of ${\CR}_{SW}(Y, \e)$
is defined by its
spectral flow which depends on the perturbation homotopy class of $\e$.
(This is different phenomenon from the (instanton) Casson invariant of
integral homology 3-spheres.)
\qed

Note that the monopole number $\# {\CR}^*_{SW}(Y, \e)$ (counted with sign)
is not a topological invariant. The number $\# {\CR}^*_{SW}(Y, \e)$
depends on the metric with harmonic-spinor jumps.

\section{Monopole homology of integral homology 3-spheres}

For an admissible perturbation $\e \in {\CP}_Y$,
we obtain a new gradient vector field
$f_{\e}$ for which the irreducibles
are all nondegenerate in \S 4. Since zeros of $f_{\e}$ are now
isolated finite-many points, we use them to generate the monopole chain groups. 

\begin{df} \label{53}
Let $(a, \p, \psi)$ and $(a^{'}, \p^{'}, \psi^{'})$ be zeros of $f_{\e}$.
A {\it chain solution}
$((A_1, \Phi_1), ... , (A_n, \Phi_n))$ from $(a, \p, \psi)$ to 
$(a^{'}, \p^{'}, \psi^{'})$ is a finite set
of Seiberg-Witten solutions over $Y\x {\R}$
which converge to $c_{i-1},
c_i \in f^{-1}_{\e}(0)$ as $t \to \mp\infty$ such that
$(a, \p, \psi) = c_0$, $c_{n} = (a^{'}, \p^{'}, \psi^{'})$,
and $(A_i, \Phi_i) \in {\CM}_{Y \x {\R}}(c_{i-1}, c_i)$ for $0 \leq i \leq n$.
\end{df}
We say that the sequence $\{(A_{\a}, \Phi_{\a})\} \in {\CM}_{Y\x {\R}}
((a, \p, \psi), (a^{'}, \p^{'}, \psi^{'}))$ is 
{\it (weakly) convergent} to the chain solution 
$((A_1, \Phi_1), ... , (A_n, \Phi_n))$ if there is a 
sequence of n-tuples of real numbers
$\{t_{\a, 1} \leq \dots \leq
t_{\a,n}\}_\a$, such that $t_{\a,i}-t_{\a,i-1}\to\infty$
as $\a \to\infty$, and
if, for each $i$, the translates
$t_{\a,i}^*(A_{\a}, \Phi_{\a}) = (A_{\a}(\circ-t_{\a,i}), 
\Phi_{\a}(\circ-t_{\a,i}))$ converge weakly to $(A_i, \Phi_i)$.

\begin{thm} \label{ufcom}
Let $\{(A_{\a}, \Phi_{\a})\} \in 
{\CM}_{Y\x {\R}}((a, \p, \psi), (a^{'}, \p^{'}, \psi^{'}))$ be a sequence of 
Seiberg-Witten solutions with uniformly bounded action over
$Y\x {\R}$. Then there exists a subsequence converging to a chain
solution $((A_1, \Phi_1), ... , (A_n, \Phi_n))$ such that
\[\mbox{Ind} D_{A_{\a}, \Phi_{\a}} = \sum_{i=1}^n
\mbox{Ind} D_{A_i, \Phi_i} = \sum_{i=1}^n (\mu_{\e}(c_i) - \mu_{\e}(c_{i-1})).\]
\end{thm}
\noindent{\bf Proof:}
It follows from the same proof as in \cite{fl} \S 3 and \cite{km},
and the compactness of Seiberg-Witten
moduli space on 4-dimensional manifolds. \qed

\begin{pro} \label{1c1}
The compactification of ${\CM}_{Y\x {\R}}(c_0, c_{n+1})$
with only chain solutions can
be described as
\[\ov{{\CM}_{Y\x {\R}}(c_0, c_{n+1})} =
\cup (\x_{i=1}^{n+1} {\CM}_{Y\x {\R}}(c_{i-1}, c_i)),\]
the union over all sequence $c_0, c_1, \cdots,
c_{n+1} \in {\CR}^*_{SW}(Y, \e)$ such that
${\CM}_{Y\x {\R}}(c_{i-1}, c_i)$ is nonempty for all
$1 \leq i \leq n+1$.

For any sequence $c_0, c_1, \cdots, c_{n+1} \in {\CR}_{SW}^*(Y, \e)$,
there is a gluing map
\[G: \x_{i=1}^{n+1} {\hM}_{Y\x {\R}}(c_{i-1}, c_i) \x
\D^{n+1} \to \ov{{\CM}_{Y\x {\R}}(c_0, c_{n+1})},\]
where $\D^{n+1}=\{(\lam_0, \cdots, \lam_n) \in [-\infty, \infty]^{n+1}:
1+\lam_{i-1} < \lam_i, 1 \leq i \leq n\}$.
\begin{enumerate}
\item The image of $G$ is a neighborhood of
$\x_{i=1}^{n+1} {\hM}_{Y\x {\R}}(c_{i-1}, c_i)$ in the
compactification with chain solutions.
\item The restriction of $G$ to
$\x_{i=1}^{n+1} {\hM}_{Y\x {\R}}(c_{i-1}, c_i) \x \mbox{Int}\,
(\D^{n+1})$ is an orientation-preserving diffeomorphism onto its image.
\end{enumerate}
\end{pro}
\noindent{Proof:} Since there is no bubbling in the Seiberg-Witten moduli space,
the map $G$ is the well-known transitivity in finite-dimensional Morse-Smale
theory. \qed

Let ${\CR}^n_{SW}(Y, \e)$ be the set of irreducible zeros $(a, \p, \psi)$
of $f_{\e}$ whose relative index $\mu_{\e} (a, \p, \psi) - 
\mu_{\e} (\T) = n$.
The {\bf monopole chain group} $MC_n(Y, \e)$
is defined to be the free Abelian
group generated by ${\CR}^n_{SW}(Y, \e)$, where the admissible 
perturbation $\e$
specifies the spectral flow $\mu_{\e}(\T)$. We write $I_{\e}(\T; \e_0)$ 
to be the
integer $\mu_{\e}(\T) - \mu_{\e_0}(\T)$ with respect to a reference 
$\e_0 \in {\CP}_Y$. Hence $\mu_{\e}(\T)$ is fixed with the fixation of
$I_{\e}(\T; \e_0)$.
Define the boundary operator $\partial : MC_n(Y, \e) \to
MC_{n-1}(Y, \e)$:
\[\partial(a, \phi, \psi) = \sum_{(a^{'}, \phi^{'}, \psi^{'}) \in 
MC_{n-1}(Y, \e)} 
\# \hat{{\cal
M}}^1_{SW, Y \x {\R}}((a, \phi, \psi), (a^{'}, \phi^{'}, \psi^{'})) \cdot 
(a^{'}, \phi^{'}, \psi^{'}). \]

\begin{pro} Let $\partial : MC_n(Y, \e) \to MC_{n-1}(Y, \e)$
be defined as above. Then
$\partial \circ \partial = 0$.
\label{5.2} \end{pro}
\noindent{\bf Proof:}
The proof follows the same argument as in (\cite{fl}, Theorem 2)
except that we have to rule out the possibility of reducible connections
entering into the picture. Note that
\[\partial^2(c_0) =  \sum_{c_1 \in 
{\CR}_{SW}^{n-1}(Y, \e)} \sum_{c_2
\in {\CR}_{SW}^{n-2}(Y, \e)}\# \hat{{\CM}}^1_{Y \x {\R}}(c_0, c_1)
\cdot  \# \hat{{\CM}}^1_{Y \x {\R}}(c_1, c_2) c_2,\]
where $c_i = (a_i, \p_i, \psi_i) \in {\CR}^*_{SW}(Y, \e) (i=0, 1, 2)$.
Consider in this sum all the terms associated to a fixed 
$c_2 \in {\CR}_{SW}^{n-2}(Y, \e)$. For the
pair $(c_0, c_2)$, 
there is the 2-dimensional moduli space ${\CM}^2_{Y \x {\R}}
(c_0, c_2)$. By Proposition~\ref{1c1},
the ends of 
$\hat{{\CM}}^2_{Y \x {\R}}(c_0, c_2)$
consists of all the components $\hat{{\CM}}^1_{Y \x {\R}}(c_0, c_1)
\x \hat{{\CM}}^1_{Y \x {\R}}(c_1, c_2)$ with
$c_1 \in {\CR}_{SW}^{n-1}(Y, \e)$. It is impossible for $c_1$ to be the
$U(1)$-reducible zero of $f_{\e}$ 
because the isotropy subgroup $\Gamma_{c_1}$ would
add to the gluing parameter and as a result would contradict the dimension
count by Proposition~\ref{sfi} and Proposition~\ref{1pf}. 
Thus
\[\sum_{c_1 \in {\CR}_{SW}^{n-1}(Y, \e)}\# \hat{{\CM}}^1_{Y \x {\R}}
(c_0, c_1) \cdot  \# \hat{{\CM}}^1_{Y \x {\R}}(c_1, c_2) =
\partial \hat{{\CM}}^2_{Y \x {\R}}(c_0, c_2) = 0. \]
\qed

As a consequence of Proposition~\ref{5.2}, for a given integral 
homology 3-sphere $Y$ and an admissible
data $\e \in {\CP}_Y$,
we have a well-defined definition of a {\bf Monopole Homology}
\[MH_*(Y; \e) = \ker \partial_* / \mbox{Im} \partial_{*+1}, \ \ \ \
* \in {\Z}.\]

Now the monopole homology $MH_*(Y; \e)$ is
sensitive to the number $I_{\e}(\T; \e_0)$, and
$MH_*(Y; \e)$ is not a topological invariant since its
Euler characteristic $\# {\CR}^*_{SW}(Y, \e)$ is metric-dependent.

\section{Homomorphisms induced by cobordisms}

From the troublesome path of metrics in $\Si_Y$
of creating/destroying harmonic
spinors (see \cite{hi}), the invariance
of the monopole homology of integral homology 3-spheres is in question.
The cobordism argument used in \cite{fl} does not apply here.
We have to construct a different cobordism
between metrics and admissible perturbations
with the fixed spectral flow $I_{\e}(\T; \e_0) 
= \mu_{\e} (\T) - \mu_{\e_0}(\T)$.
In this section, we show that
our monopole homology is independent of metrics and of admissible
perturbations within the class $I_{\e}(\T; \e_0)$.

Let $X$ be an oriented 4-manifold with two cylindrical ends $Y_1 \times {\bf
 R}_+$ and $Y_2 \times {\bf R}_-$, where $Y_1$ and $Y_2$
are integral homology 3-spheres.
Let $\tau : X \to [0, \infty )$ be a smooth cutoff
function such that $\tau(x) = 0$ for $x$ lying outside of $Y_1 \times {\bf
R}_+ \cup Y_2 \times {\bf R}_-$ and $\tau (y,t) = | t|$ for $(y,t) \in Y_1
\times {\bf R}_+ \cup Y_2 \times {\bf R}_-$ and $|t| > t_0 > 0$ and
$e_{\delta } = e^{\delta \tau (x)}$. Then
using the cutoff function $\tau $ and a background connection  we can
extend $\frac{d}{dt} + \alpha , \frac{d}{dt} +
\beta $ to a connection $\nabla_0$ on $X$ such that
 \[\nabla_0|_{Y_1 \times [t_0 , \infty)} = \frac{d}{dt} +\alpha , \ \ \
 \nabla_0|_{Y_2\times (- \infty, - t_0]} =\frac{d}{dt} +
 \beta . \] Similarly, we can extend sections on $W_X^{\pm}$.
The Fr\'{e}chet space $\Omega^1_{\mbox{comp}}(X, AdP)\op 
\G_{\mbox{comp}}(W_X^{\pm})$ of compact supported
$C^{\infty }$-sections on $(T^*X \otimes AdP)\op \G (W_X^{\pm})$ 
can be completed to a Banach space
\[{\cal A}^p_{k, \delta} (X) = (\nabla_0, 0) + L^p_{k, \delta }
(\Omega^1(X, AdP)\op \G (W_X^{\pm})),\]
where $\|c\|_{L^p_{k, \delta }} = \|e_{\delta } \cdot c \|_{
L^p_k}$ for 
$c \in \Omega^1_{\mbox{comp}}(X, AdP)\op \G_{\mbox{comp}}(W_X^{\pm})$. 
The gauge group ${\cal G}^p_{k+1, \delta}$
is given by
$L^p_{k+1, \delta}$-norm of $\mbox{Aut}(\det W_X^{\pm})$.
So the quotient space is 
${\cal B}^p_{k, \delta }(X) = {\cal A}^p_{k, \delta} (X)/
{\cal G}^p_{k+1, \delta}$.
The perturbation data $\e_1 = (g_{Y_1}, \a_1)$ and $\e_2 = (g_{Y_2}, \a_2)$
at the ends provide 
the gradient vector fields $f_{\e_1}$ and $f_{\e_2}$ so
that the zeros of $f_{\e_1}$ on $Y_1$ and of $f_{\e_2}$ on $Y_2$ are generic.
Clearly these
perturbation data $\e_1$ and $\e_2$ 
can be pulled back to the cylindrical
ends $Y_1 \times {\bf R}_+$ and
$Y_2 \times {\bf R}_-$, and produce perturbations on the
time-invariant monopole equation on ${\cal B}^p_{k, \delta}(Y_1
\times {\bf R}_+ )$ and ${\cal B}^p_{k, \delta}(Y_2 \times {\bf R}_-)$ (same
$\delta $ as before).
According to (\cite{fl} (1c.2) and \cite{km, tau, wi}), 
there exists a Baire's first category
subset in the space ${\cal M}et(X) \times \Pi_X$ of Riemannian metrics $g_X$ and
perturbation data $\a_X$ such that ${\cal M}_{\e_X}(c, c^{'})$ 
($\e_X = (g_X, \a_X)$)
is a smooth manifold with 
\begin{equation} \label{1dim}
\mbox{dim} {\cal M}_{\e_X}(c, c^{'}) =\mu_{\e_1} (c) - \mu_{\e_2} (c^{'}) 
+ \frac{1}{2}(2 \c + 3\s)(X). 
\end{equation}
In addition, ${\cal M}_{\e_X}(c, c^{'})$
is oriented with an
orientation specified by the orientations on
$H^1(X, \bf R)$ and $H^0(X, {\bf R}) \oplus H^2_+(X, \bf R)$ 
(see \cite{don, km, tau, wi}).

Define a homomorphism $\Psi_* = \Psi_*(X; \e_X): MC_*(Y_1; 
\e_1) \to MC_*(Y_2; \e_2)$
of the monopole chain complexes by the
formula
\[\Psi_*(c) = \sum_{c^{'} \in {\cal R}^*_{SW}(Y_2, \e_2)} 
\# {\cal M}^0_{\e_X} (c, c^{'}) \cdot c^{'}, \ \ \ \ 
c \in {\cal R}^*_{SW}(Y_1, \e_1),\]
where ${\cal M}^0_{\e_X} (c, c^{'})$ is the $0$-dimensional
oriented moduli space connecting $c$ to $c^{'}$ on $X$ and
$\mu_{\e_1}(c) -
\mu_{\e_2}(c^{'}) = - \frac{1}{2}(2 \c + 3\s)(X)$.

\begin{pro} \label{corb}
Given a cobordism $X$ and perturbation data $\e_X
\in {\cal M}et(X) \times \Pi_X$ as before, the homomorphism $\Psi_*$ is a
chain map shifting the degree by $\frac{1}{2}(2 \c + 3\s)(X)$. Furthermore the
induced homomorphism \[\Psi_* = \Psi_* (X; \e_X): MH_*(Y_1;
\e_1) \to MH_*(Y_2; \e_2) \]
on the monopole homologies depends only on the cobordism $X$.
\end{pro}
\noindent{\bf Proof:} It follows the same argument as in 
\cite{fl} Theorem 3 and \cite{ll} \S 5. \qed
\vspace{0.1in}

We show below that $\Psi_*(X; \e_X)$ is functorial with respect to
the composite cobordism. Given two cobordisms $(U; \e_U)$ connecting
$Y_1$ to $Y_2$ and $(V; \e_V)$ connecting $Y_2$ to $Y_3$ so that
$\e_U$ and $\e_V$ agree on $Y_2$, we can form the
composite cobordism $(W; \e_W)$ connecting $Y_1$ to $Y_3$. Then
\begin{equation} \label{compo}
\Psi_*(W; \e_W) = \Psi_*(V; \e_V) \circ
\Psi_*(U; \e_U).
\end{equation}
A different strategy from Floer's has to be taken to prove that $MH_*(Y,
\e)$ is independent of admissible perturbations 
$\e = (g_Y, \a)$ within the class of $I_{\e}(\T; \e_0)$. 
We consider the
time-dependent perturbations of the Seiberg-Witten equation and
its associated moduli space. 
Given two admissible perturbation data
of generic metrics $g^{-1}_Y$ and $g_Y^1$ and 1-forms 
$\a_{-1}$ and $\a_1$ with $I_{\e_{-1}}(\T; \e_0) = I_{\e_1}(\T; \e_0)$
(here $\e_t = (g_Y^t, \a_t)$),
there is an one-parameter family of admissible perturbations $\Lambda =
\{\e_t = (g_Y^t, \a_t)| - \infty \leq t \leq \infty\}$ joining them.
Assume that the pair $\e_t = (g_Y^{-1}, \a_{-1})$
for $t \leq -1$ and $\e_t = (g_Y^1, \a_1)$
for $t \geq 1$. On the cylinder $Y \times \bf R$, we consider the
perturbed Seiberg-Witten equation
\begin{equation} \label{ASDp}
\frac{\bd \psi}{\bd t} + \bd_{a_t}^{\nabla_{g_Y^t}+\a_t} \psi = 0, \ \ \ 
\frac{\partial a_t}{\partial t} + *_{g_Y^t} F(a_t) + \a_t =
i\tau_{g_Y^t}(\psi, \psi). \end{equation}
Given $c \in {\cal R}^*_{SW}(Y, \e_{-1})$ and
$c^{'} \in {\cal R}^*_{SW}(Y, \e_1)$,
we denote by ${\cal M}_{\Lambda }(c, c^{'})$ the subspace in ${\cal
B}^p_{k,\delta }(c, c^{'})$ consisting of solutions of (\ref{ASDp}).
Then there exists a homomorphism
\[ \Psi_{\Lambda }: MC_n (Y; \e_{-1}) \to MC_n (Y; \e_1) \]
of the monopole chain complexes defined by
\[ \Psi_{\Lambda }(c) = \sum_{c^{'} \in {\cal R}^n_{SW}(Y, \e_1)} \# {\cal
M}_{\Lambda }^0(c, c^{'}) \cdot c^{'}, \ \ \ \ c \in {\cal
R}^n_{SW}(Y, \e_{-1}). \]

\begin{pro} \label{6.5}
Let $\Lambda =
\{\e_t = (g_Y^t, \a_t)| t \in {\bf R}\}$ be an
family of admissible perturbations as defined above such that
$Ind D_{\e_t}(\T) = 0$. Then
\begin{enumerate}
\item If $\Lambda $ is a constant family of admissible perturbations
($g_Y^t = g_Y, \a_t = \a$), then $\Psi_{\Lambda } = id$.
\item $\Psi_{\Lambda }$ is a chain map: $\partial \Psi_{\Lambda } =
\Psi_{\Lambda }\partial $.
\item Given two families $\Lambda$ and $\Lambda^{'}$ of
admissible perturbations
joining $(g_Y^{-1}, \a_{-1})$ to $(g_Y^0, \a_0)$ and from
$(g_Y^0, \a_0)$ to $(g_Y^1, \a_1)$, we have
$\Psi_{\Lambda \circ \Lambda^{'}} =
\Psi_{\Lambda } \circ \Psi_{\Lambda^{'}}.$ 
\item If a family $\Lambda_0$ of admissible perturbations connecting
$(g_Y^{-1}, \a_{-1})$ and $(g_Y^1, \a_1)$ can be deformed into another
$\Lambda_1$ by admissible families $\Lambda_{\lambda} (0\leq \lambda \leq
1)$, then the two monopole 
chain maps $\Psi_{\Lambda_0}$ and $\Psi_{\Lambda_1}$ are chain
homotopic to each other.
\end{enumerate}
\end{pro}
\noindent{\bf Proof:} (1) If the perturbation is
time independent $\e_t = (g_Y, \a)$, then ${\cal
M}_{\Lambda }^0(c, c^{'})$ is just the space ${\cal M}_{Y \x {\bf R}}^0
(c, c^{'})$. 
For the 0-dimensional component ${\cal
M}^0_{\Lambda }(c, c^{'})$, this means time-invariant solutions $c_t$
on $Y \x \bf R$, and 
we have $[c_t] = c = c^{'}$. Therefore
$\#{\cal M}^0_{\Lambda }(c, c^{'}) =
\delta_{c c^{'}}$ and $\Psi_{\Lambda} = id$.

(2) We consider the compactification of ${\cal
M}_{\Lambda }(c, c^{'})$ as developed in \cite{fu, ll}.
By Proposition~\ref{1c1} and \cite{km, tau, wi}, ${\cal
M}_{\Lambda }(\alpha , \beta)$ can be compactified such that the
codimension-one boundary consists of
\begin{equation}
\cup_{c_{-1}} \hat{{\cal M}}_{Y \x {\bf R}}(c, c_{-1}) \times_{c_{-1}}
{\cal M}_{\Lambda
}(c_{-1}, c^{'}) \coprod \cup_{c_1}{\cal M}_{\Lambda }(c, c_1)
\times_{c_1} \hat{{\cal M}}_{Y \x {\bf R}}(c_1, c^{'}) .
\label{6.6} \end{equation}
Here $c_{\pm 1} \in {\cal R}_{SW}(Y, \e_{\pm 1})$ and 
${\cal M}_{Y \x {\bf R}}(c, c_{-1})$ is
the moduli space of monopoles on $Y \times (- \infty , -1)$ with
respect to the perturbation $\e_{-1}$ and $\hat{{\cal
M}}_{Y \x {\bf R}}(c, c_{-1}) = {\cal M}_{Y \x {\bf R}}(c, c_{-1})/{\bf R}$.
Similarly $\hat{{\cal
M}}_{Y \x {\bf R}}(c_1, c^{'})$ is obtained from the perturbation data $\e_1$.
Consider
the 1-dimensional components ${\cal M}^1_{\Lambda }(c, c^{'})$ of ${\cal
M}_{\Lambda }(c, c^{'})$, whose boundary by (\ref{6.6}) gives two
types of oriented points counted as $\bd \Psi_{\Lambda } =  \Psi_{\Lambda }\bd$.
We can rule out the possibilities of the reducible $\T$
for $c_{\pm 1}$. If they occurred, 
then they would have an
additional $U(1)$-symmetry on these moduli spaces. This is impossible by the
dimension reasoning from Proposition~\ref{sfi}, Proposition~\ref{1pf} and
our hypothesis $I_{\e_{-1}}(\T; \e_0) = I_{\e_1}(\T; \e_0)$ (see below also).

(3) For a composite cobordism and its induced
homomorphism, we study the moduli space ${\cal
M}_{\Lambda * \Lambda^{'}}(T; \alpha , \beta)$ of solutions of the 
Seiberg-Witten
equation on $Y \x {\bf R}$ with respect to the following time-dependent
admissible perturbation
data $\Lambda *_T \Lambda^{'}$, where
\[\Lambda *_T \Lambda^{'} = \left\{ \begin{array}{lc}
\e_{-1} = (g_Y^{-1}, \a_{-1}) & - \infty < t \leq -T -1\\
\Lambda = (g_Y^{t+T}, \a_{t+T}) & -T -1 \leq t \leq -T\\
\e_0 & -T \leq t \leq T\\
\Lambda^{'} = (g_Y^{t-T}, \a_{t-T})& T \leq t \leq T+1\\
\e_1 & T+1 \leq t < + \infty.  
\end{array} \right. \]
Let $T$ be
sufficiently large. Thus ${\cal M}_{\Lambda * \Lambda^{'}}(T; c, c^{'}) 
(T \geq T_0)$ is approximated by the union
\begin{equation}
\cup_{c_0} \overline{{\cal M}}_{\Lambda }(c, c_0)
\times_{c_0}\overline{{\cal M}}_{\Lambda^{'}}(c_0, c^{'}).
\label{6.7} \end{equation}
where $\overline{{\cal M}}_{\Lambda }(c, c_0)={\cal M}_{\Lambda
}(c, c_0)/(\Gamma_c \x \Gamma_{c_0})$. Note that
the 0-dimensional components in $\overline{{\cal M}}_{\Lambda }(c, c_0)
\times_{c_0}\overline{{\cal M}}_{\Lambda^{'}}(c_0, c^{'})$
correspond to the $c^{'}$-coefficients in
\[\Psi_{\Lambda^{'}} \circ \Psi_{\Lambda }(c) = \sum_{c_0} \#
\overline{{\cal M}}_{\Lambda }^0(c, c_0)
\cdot \#\overline{{\cal M}}^0_{\Lambda^{'}}(c_0, c^{'}) \cdot c^{'}.\]
On the other hand, as $T \to 0$, the 0-dimensional
component of the moduli space ${\cal M}_{\Lambda *
\Lambda^{'}}(T; c, c^{'})$
gives the $c^{'}$-coefficients
in $\Psi_{\Lambda * \Lambda^{'}} (c) = \sum
{\cal M}^0_{\Lambda * \Lambda^{'}}(c, c^{'})\cdot c^{'} $. Because
$\cup_{0\leq T \leq T_0} {\cal M}^0_{\Lambda *
\Lambda^{'}}(T; c, c^{'})$ is the cobordism between
${\cal M}^0_{\Lambda *\Lambda^{'}}(0; c, c^{'})$ and ${\cal
M}^0_{\Lambda *\Lambda^{'}}(T_0; c, c^{'})$, so the assertion (3) follows 
by ruling out the reducible $\T$. Note that
\[\mbox{dim}\overline{{\cal M}}_{\Lambda} (c, c_0) =
\mu_{\e_{-1}}(c) -
\lim_{\e_t \in {\Lambda }, \e_t \to \e_0} \mu_{\e_t}(c_0)
 - \mbox{dim}\Gamma_{c_0}; \]
\begin{equation} \label{dim}
\mbox{dim}\overline{{\cal
M}}_{\Lambda^{'}} (c_0, c^{'}) =
\lim_{\e_t \in {\Lambda^{'}}, \e_t \to \e_0}
\mu_{\e_t }(c_0) - \mu_{\e_1}(c^{'}).
\end{equation}
By Proposition~\ref{sfi} and Proposition~\ref{1pf}, we obtain 
\[\lim_{\e_t \in {\Lambda }, \e_t \to \e_0} \mu_{\e_t}(c_0) =
\lim_{\e_t \in {\Lambda^{'}}, \e_t \to \e_0}
\mu_{\e_t}(c_0) = \mu (c_0).\]
So it satisfies the
equations $\mu_{\e_{-1}}(c) -
\mu (c_0) = 1$ ($c_0 = \T$) and $\mu (c_0) -
\mu_{\e_1}(c^{'}) = 0$.
This is impossible because of $\mu_{\e_{-1}}
(c) = \mu_{\e_1}(c^{'})$.
{\em If these spectral flows $I_{\e_{\pm 1}}(\T; \e_0)$ 
are not fixed to be same, then the above argument becomes
invalid. }

(4) Let $\Lambda_i (i=0, 1)$ be a family of 
time-independent admissible perturbations which
connect up $\e_{-1}$ and $\e_1$.
Suppose that $\Lambda_0$ and $\Lambda_1$
can be smoothly deformed from one to another by a 1-parameter family
$\Lambda_s= \{\e_t^s = (g_Y^{s,t}, \a^s_t) , 0\leq s\leq 1, \ -1\leq t \leq 1\}$
of the same type of admissible perturbations. Set $\Lambda_s =
\Lambda_0$ for $0 \leq s \leq \frac{1}{4}$ and $\Lambda_s = \Lambda_1$ for
$\frac{3}{4} \leq s \leq 1$. Associated to this situation, there is a
1-parameter family of moduli spaces denoted by
${\cal H}\tilde{{\cal M}}(c, c^{'})
= \cup_{0\leq s \leq 1}\tilde{{\cal M}}_{\Lambda_s}(c, c^{'})$,
\[{\cal H}\tilde{{\cal M}}(c, c^{'}) = \{(\Phi, s) | \Phi \in
\tilde{{\cal M}}_{\Lambda_s}(c, c^{'}), 0\leq s \leq 1 \} \subset
{\cal B}_{k, \d}^p(c, c^{'}) \times [0,1], \]
where ${\cal H}\tilde{{\cal M}}$ is the set of regular solutions of
Seiberg-Witten equation with respect to $\e_t^s$,
and is a smooth manifold with dimension
$\mu_{\e_{-1}}(c) -
\mu_{\e_{1}}(c^{'}) +1$.
The codimension-one boundary consists of
\[{\cal M}_{\Lambda_1}(c, c^{'}) \times \{0\} \coprod
{\cal M}_{\Lambda_0}(c, c^{'}) \times \{1\}, \]
\[\cup_{(s, c_0)}\tilde{{\cal M}}_{\Lambda_s}(c, c_0)
\times {\cal M}_{\e_1}(c_0, c^{'}) \coprod
\cup_{(s, \gamma )}{{\cal
M}}_{\e_{-1}}(c, c_0)
\times \tilde{{\cal M}}_{\Lambda_s}(c_0, c^{'}) .\]
Since $\tilde{{\cal M}}_{\Lambda_s}(c, c_0)$ and $\tilde{{\cal
M}}_{\Lambda_s}(c_0, c^{'})$
are solutions of the Seiberg-Witten equation with
virtual dimension $-1$, they can only occur for $0<s<1$.
The homomorphism $H: MC_*(Y; \e_{-1})
\to MC_*(Y; \e_1)$ of degree $+1$ is defined by
\[H(c) = \sum_{c_0} \sum_s \# \tilde{{\cal
M}}^0_{\Lambda_s}(c, c_0)\cdot c_0, \ \ \ \mbox{for} \ \
c \in {\cal R}^n_{SW}(Y, \e_{-1}), c_0 \in {\cal R}^{n+1}_{SW}(Y, \e_1). \]
That $c_0$ is reducible is
eliminated by the extra $U(1)$-symmetries in ${\cal M}_{\e_1}(c_0, c^{'})$
and ${\cal M}_{\e_{-1}}(c, c_0)$ and 
$I_{\e_1}(\T; \e_0) = I_{\e_{-1}}(\T; \e_0)$.
Summing up $c^{'} \in {\cal R}^n_{SW}(Y, \e_1)$, we have
\[ \Psi_{\Lambda_0} (c) - \Psi_{\Lambda_1}(c) = H\circ
\partial_{\e_{-1}}(c) + \partial_{\e_1} \circ H(c).\]
Therefore $\Psi_{\Lambda_0}$ and $\Psi_{\Lambda_1}$ are monopole
chain homotopic to each other.
\qed

Thus the monopole 
homology groups $MH_*(Y; \e^{\pm 1})$
associated to two admissible perturbation data
are canonically isomorphic to each other
whenever $I_{\e^1}(\T; \e_0) =
I_{\e^{-1}}(\T; \e_0)$ for the unique $U(1)$-reducible $\Theta$ on $Y$.
Thus it is more appropriate
to denote $MH_*(Y; \e)$ by $MH_*(Y; I_{\e}(\T; \e_0))$. 
For an integral homology 3-sphere $Y$, the monopole
homology can be extended to a function $$MH_{SWF}: \{I_{\e}(\T; \e_0): \e \in
{\CP}_Y\} \to \{MH_*(Y, I_{\e}(\T; \e_0)): \e \in {\CP}_Y\}.$$
(Changing a reference $\e_0$ corresponds to the same homology groups
with grading $I_{\e_0^{'}}(\T; \e_0)$-shift)
This function $MH_{SWF}$ is a topological invariant of the
integral homology 3-sphere $Y$, up to the degree-shifting of monopole
homologies. Hence
such a function $MH_{SWF}$ may be called a Seiberg-Witten-Floer theory,
which is completely different from the
instanton Floer homology, but more related to the treatment in \cite{ll}.

\section{Relative Seiberg-Witten invariants}

The Seiberg-Witten invariant (see
\cite{don, tau, wi}) has proved so useful and at least powerful as the
Donaldson
invariant in many cases, and is much easier to compute.
In this section we are going to extend the Seiberg-Witten invariant to the
relative one on 
smooth 4-manifolds with boundary integral homology 3-spheres. The ``relative
Seiberg-Witten invariants'' is no longer a topological invariant since it lies 
in a monopole homology depending upon Riemannian metrics of integral
homology 3-spheres. But the natural pairing between
``relative Seiberg-Witten invariants'' does recover the 
Seiberg-Witten invariant 
of closed smooth 4-manifolds.

Let $X$ be a smooth 4-manifold with $b_1(X) > 0$ and
boundary $Y$ (an integral homology
3-sphere). The collar of $X$ can be identified with $Y \times [-1,1]$, and
the admissible perturbation data on $Y$ can be extended
inside $X$ as we did in \S 7.
Fixing $I_{\e}(\T; \e_0)$ should be understood though this section.

\begin{df} For a smooth 4-manifold $X$ with boundary $Y$ (an integral homology
3-sphere), the 0-degree relative Seiberg-Witten invariant is defined by
\[q_{X,Y, \e} = \sum_{c \in {\cal R}^*_{SW}(Y, \e)} \# {\cal M}^0_X(c)\cdot
c, \]
where ${\cal R}^*_{SW}(Y, \e)$ is the set of all nondegenerate zeros of
$f_{\e}$ with 
prescribed $I_{\e}(\T; \e_0)$.
\end{df}  \label{relative}

By the index calculation and our convention $\mu_{\e}(c) = 
SF(c, \T)$, we have
\[\mbox{dim}{\cal M}_X^0(c) + \mu_{\e}(c) = \mbox{dim} {\cal M}_X(\T) = 
\frac{1}{4}(c_1(\pi^*(L))^2 - (2 \c+ 3\s))(X) =
- \frac{1}{4}(2 \c+ 3\s)(X),\]
since $c_1(L) =0$ for the integral homology 3-sphere $Y$.
Thus $q_{X,Y, \e}$ is in the monopole chain group with grading
$- \frac{1}{4}(2 \c+ 3\s)(X)$.

\begin{pro} \label{cycle}
For $q_{X,Y, \e} \in MC_{\mu_X}(Y, \e)$ with 
$\mu_X = - \frac{1}{4}(2 \c+ 3\s)(X)$ and a fixed class $I_{\e}(\T; \e_0)$,
we have $\bd_Y \circ q_{X,Y, \e} = 0$.
\end{pro}
\noindent{\bf Proof:} \[\partial_Y \circ q_{X,Y, \e}(c) = \sum_{c \in 
{\cal R}^{\mu}_{SW}(Y, \e)}
\sum_{c^{'} \in {\cal R}^{\mu -1}_{SW}(Y, \e)} \# {\cal
M}_X^0(c) \cdot \# \hat{{\cal M}}^1_{Y \times {\bf R}}(c, c^{'})
\cdot c^{'}. \]
For both $c$ and $c^{'}$ irreducible (nondegenerate) zeros of $f_{\e}$, 
we take one-dimensional moduli space ${\cal M}^1_X(c^{'})$ for
fixed $c^{'}$. Then we count the ends of the moduli space to conclude the
result. Again it is a technical point to avoid the reducible $\T$
entering the boundary ${\cal M}_X(\T) \times
{\cal M}_{Y \times {\bf R}}(\T, c^{'})$. For
the reducible $\T$, we have the dimension counting 
\[ \mbox{dim} \{{\cal M}_X(\T) \x {\cal M}_{Y \times {\bf
R}}(\T, c^{'})\} = \mbox{dim} {\cal M}_X(\T) + \mbox{dim}\G_{\T} +
\mbox{dim}{\cal M}_{Y \times {\bf R}}(\T, c^{'}) \geq 0+1+1 = 2. \]
So $c$ cannot be the reducible $\T$, and $\partial_Y \circ q_{X,Y, \e} = 0$. 
Hence $q_{X,Y, \e}$ is indeed a
monopole cycle. \qed

Let $q_{X,Y, \e}(g_X)$ be the relative Seiberg-Witten invariant with respect
to the metric $g_X$. Now we show that the monopole homology class
$[q_{X,Y, \e}(g_X)]$ defined by Proposition~\ref{cycle}
is independent of
metrics $g_X$ with $g_X|_Y$ in the fixed class of $I_{\e}(\T; \e_0)$.

\begin{pro} \label{metri}
Let $g_X^i (i = 1, 2)$ be two generic metrics on $X$ with induced
metric $g_Y^i$ generic such that $I_{\e_1}(\T; \e_0) = I_{\e_2}(\T; \e_0)$ and 
$\e_i = (g_Y^i, \a_i)$. Then there exist $c^{'} \in MC_{\mu_X + 1}$ with 
$\mu_X = - \frac{1}{4}(2 \c+ 3\s)(X)$
such that we have
\[ q_{X,Y, \e_2}(g_X^2) - q_{X,Y, \e_1}(g_X^1) = \partial (c^{'}). \]
In particular, $[q_{X,Y,\e_2}(g_X^2)] = [q_{X,Y, \e_1}(g_X^1)]$
as the monopole 
homology class in $MH_{\mu_X}(Y, I_{\e_i}(\T; \e_0))$.
\end{pro}
\noindent{\bf Proof:}
Let $\{g_X^{t+1}\}_{0 \leq t \leq 1}$ be a family of metrics on $X$
such that $I_{\e_{t+1}}(\T; \e_0)$ is independent of $t$ with
$\e_{t+1} = (g_X^{t+1}|_Y, \a_{t+1})$ and ${\cal M}_X^0(g_X^{t+1})(c)$
has virtual dimension 0 with respect to $c$
irreducible. Therefore
$\{{\cal M}_X^0(g_X^{t+1})(c)\}_{0 \leq t \leq 1}$
is an one-dimensional moduli space of Seiberg-Witten solutions on
$X$. The corresponding codimension-one
boundary in $[0,1]
\times {\cal B}_X(g_X^{t+1})(c)$ is given by
\[\partial (\{{\cal M}_X^0(g_X^{t+1})(c)\}_{0 \leq t \leq 1})= \]
\[\{0\} \times {\cal M}_X^0(g_X^1)(c) \coprod - \{1\} \times {\cal
M}_X^0(g_X^2)(c) \coprod \partial( \sum_{\mu_{\e_{t+1}}(c) -
\mu_{\e_{t+1}}(c^{'}) = - 1}
\# ([0,1] \x {\cal M}^{-1}_X(g_X^{t+1})(c^{'}))). \]
The number $\langle \partial_Y c^{'}, c \rangle$ is the algebraic number
of $([0,1] \times {\cal M}^{-1}_X(g_X^{t+1})(c^{'}))$.
The $c^{'}$ cannot be the
reducible $\T$ by the fixed $I_{\e_1}(\T; \e_0)$ 
with the same argument as before.
So 
\[q_{X,Y,\e_2}(g_X^2)(c) - q_{X,Y,\e_1}(g_X^1)(c)  = \langle \partial_Y c^{'},
c \rangle.\]
Hence $q_{X,Y,\e_i}(g_X^i) (i= 1, 2)$ (as a monopole cycle) gives the same
monopole homology class. \qed

Note that orientation reversing from $Y$ to $-Y$ changes the grading from
$\mu_{\e}(c)$ to $-1 - \mu_{\e}(c)$
(certainly does not change the
solutions of the Seiberg-Witten equation on the 3-manifold), 
so there is a nature identification between
$MC_{\mu_{\e}}(Y, \e)$ and $CF_{-1- \mu_{\e}}(-Y, \e)$.

\begin{thm} \label{invariant} For a smooth 4-manifold $X =X_0 \#_Y X_1$
with $b_2^+(X_i) > 0 (i = 0, 1)$ and $Y$ an integral homology 3-sphere,
the Seiberg-Witten invariant of the 4-manifold $X$ is
given by the Kronecker pairing of $MH_*(Y; I_{\e}(\T; \e_0))$ with 
$MH_{-1-*}(-Y;I_{\e}(\T; \e_0))$ for
$q_{X_0,Y, \e}$ and $q_{X_1,- Y, \e}$;
\[ \langle  , \rangle :
MH_*(Y;I_{\e}(\T; \e_0)) \x MH_{-1-*}(-Y;I_{\e}(\T; \e_0)) \to {\bf Z}; 
\ \ \  q_{SW}(X) =
\langle q_{X_0,Y, \e}, q_{X_1,-Y, \e}\rangle. \]
More precisely,
$q_{SW}(X_0 \#_Y X_1) = \sum_{c}\# {\cal M}_{X_0,Y, \e}^0(c) \cdot \# {\cal
M}_{X_1, -Y}^0(- c)$, where $I_{\e}(\T; \e_0)$ is fixed.
The invariant $q_{SW}(X)$ is independent of the choice of $I_{\e}(\T; \e_0)$.
\end{thm}
\noindent{\bf Proof:} If $Y$ admits a metric of positive scalar
curvature, then the proof
is given in \cite{wi} with $I_{\e}(\T; \e_0) = 0$ the special case.
The assumption implies that
$b^+_2(X) > 1$. So we can rule out the existence of reducible
solutions on $X$ by the standard method (see \cite{don, km, tau, wi}). 
Note that
\[\mbox{dim} {\cal M}_{X_0}(c) + \mbox{dim} {\cal M}_{X_1}(c) + \mbox{dim}
\G_{\T}
= \mbox{dim} {\cal M}_{X}.\]
By the dimension equation,
we can eliminate the term $\# {\cal M}_{X_0,Y, \e}^0(c) \cdot
\# {\cal M}_{X_1, -Y, \e}^0(- c)$ with $c= \T$. Then the 0-dimensional
moduli space on $X$ is obtained by gluing the solutions on $(X_0, Y)$ with
ones on $(X_1, -Y)$. Using the standard technique on stretching the neck 
\cite{dk}, one gets the equality
$q_{SW}(X) = \langle q_{X_0,Y, \e}, q_{X_1,-Y, \e}\rangle$.
Since $q_{SW}(X)$ is a
topological invariant, so the pairing is independent of the choice of
$I_{\e}(\T; \e_0)$.
\qed

For higher degree relative Seiberg-Witten invariants, one can obtain
the similar results as in \cite{ll}.

Computing the monopole homology is extremely complicated due to the
Riemannian metric, harmonic spinor, spectral flow and solution of the
first-order Dirac-type nonlinear differential equation. Even for
the 3-sphere, a complete calculation of the function $MH_{SWF}$ is very 
difficult at this moment. Understand the harmonic spinors on $S^3$ with a
subfamily of Riemannian metrics (metrics are $SU(2)$-left invariant and 
$U(1)$-right invariant) is already quite involved by the work of Hitchin
\cite{hi}. On the other hand, Theorem~\ref{invariant} gives us a flexibility
to understand the Seiberg-Witten invariant of closed smooth
4-manifolds through the
relative ones with some preferred Riemannian metric(s) on the integral
homology 3-sphere.

\noindent{\bf Remark:} The method we developed in this paper also can be 
extended to rational homology 3-spheres with fixed spectral flows along all
$U(1)$-reducible solutions of Seiberg-Witten equation on the
rational homology 3-sphere (see \cite{ll} for more detail).

\noindent{\bf Acknowledgement}: The author would like to thank R. Lee for many 
discussions in our joint paper \cite{ll} which is a root for this paper. 
Realizing the correction term by the spectral flow is initiated from 
\cite{clm, ll}. 
It is a pleasure to thank Cappell, Lee and Miller whose work in \cite{clm}
inspired the circle of ideas in this paper.

\end{document}